\documentclass[draft]{article}

\usepackage[leqno]{amsmath}
\usepackage{amssymb}
\usepackage{amsthm}
\usepackage{afterpage}
\usepackage{graphics}
\theoremstyle{definition}

\newtheorem{tm}{Theorem}[section]
\newtheorem{lm}[tm]{Lemma}
\newtheorem{cor}[tm]{Corollary}

\newcommand{\fej}[1]{\section{\!\!\!\!\!\!.\ #1}}
\newcommand{\alfej}[1]{\subsection{\!\!\!#1}}
\newcommand{\nemfej}[1]{\renewcommand\thesection{}\section{\!\!\!\!\!\!#1}\renewcommand\thesection{\arabic{section}}}

\newcommand*{\hop}{\bigskip\noindent}
\newcommand*{\Zb}{\mathbb Z}
\newcommand*{\Rb}{\mathbb R}
\newcommand*{\om}{\omega}
\newcommand*{\omi}{\un\om^{(i,\,i+1)}}
\newcommand*{\omir}{\un\om^{(i+1,\,i)}}
\newcommand*{\omipp}{\un\om^{(i+1,\,+)}}
\newcommand*{\omipm}{\un\om^{(i+1,\,-)}}
\newcommand*{\omim}{\un\om^{(i,\,-)}}
\newcommand*{\omip}{\un\om^{(i,\,+)}}
\newcommand*{\omimm}{\un\om^{(i-1,\,-)}}
\newcommand*{\un}[1]{\underline{#1}}
\newcommand*{\be}{\beta}
\newcommand*{\e}[1]{\text{\rm e}^{#1}}
\newcommand*{\ve}{\varepsilon}
\newcommand*{\di}{\,\text{\rm d}}
\newcommand*{\Ev}{{\bf E}}
\newcommand*{\Evt}{{\bf E}^{(\un\te)}}
\newcommand*{\vp}{\varphi}
\newcommand*{\pt}{\partial}
\newcommand*{\te}{\theta}
\newcommand*{\Evtipp}{\Ev^{(\un\te^{(i+1,\,+)})}}
\newcommand*{\teim}{\un\te^{(i,\,-)}}
\newcommand*{\teip}{\un\te^{(i,\,+)}}
\newcommand*{\teipp}{\un\te^{(i+1,\,+)}}
\newcommand*{\teimm}{\un\te^{(i-1,\,-)}}
\newcommand*{\Evtimm}{\Ev^{(\un\te^{(i-1,\,-)})}}

\begin{document}
\title{Multiple shocks in bricklayers' model}
\author{M\'arton Bal\'azs\\
\\
Department of Mathematics,\\
University of Wisconsin-Madison}
\date{March 12,\ 2004}
\maketitle

\begin{abstract}
In bricklayers' model, which is a generalization of the misanthrope processes, we show that a nontrivial class of product distributions is closed under the time-evolution of the process. This class also includes measures fitting to shock data of the limiting PDE. In particular, we show that shocks of this type with discontinuity of size one perform ordinary nearest neighbor random walks only interacting, in an attractive way, via their jump rates. Our results are related to those of Belitsky and Sch\"utz \cite{qse} on the simple exclusion process, although we do not use quantum formalism as they do. The structures we find are described from a fixed position. Similar ones were found in Bal\'azs \cite{valak}, as seen from the random position of the second class particle.
\end{abstract}

\hop


\nemfej{Introduction}

It has been known for a while that many models from the field of stochastic interacting particle systems have deterministic partial differential equations as scaling limit, see e.g.\ the book of Kipnis and Landim \cite{cl}, or Sepp\"al\"ainen \cite{hkl}, T\'oth  and Valk\'o \cite{hydro}. In the so-called \emph{Eulerian scaling} this equation is usually a Burger-type one, developing shocks as entropy solutions from any decreasing initial data. The microscopic structure of these shocks is of great interest. Many results construct such structures from the viewpoint of the so-called \emph{second class particle}, which is an object coming from \emph{coupling} of these models. Some examples are De Masi, Kipnis, Presutti and Saada \cite{mkps}, Derrida, Lebowitz and Speer \cite{dls}, Ferrari \cite{shock}, and Ferrari, Fontes and Kohayakawa \cite{imes}. As this particle has a complex random motion in the system, it is not immediately clear how these structures look like from a fixed, non-moving position. On the other hand, this way of looking at the problem makes it difficult to deal with the case of multiple shocks. Ferrari, Fontes and Vares develop these methods and handle this case in \cite{semuso}.

Another work in this area is by Bal\'azs \cite{valak}, who introduces a simple product measure which shows the properties of a shock, and is stationary as seen by the second class particle. This result is valid for the \emph{exponential bricklayers' process} (see below for a definition). Our aim here is to identify a structure corresponding to this distribution, but this time as seen from a fixed, non moving position. This way of looking at the problem makes it easy to handle multiple shocks as well, thus we can also investigate how these shocks interact with each other at the microscopic level.

We consider the class of \emph{bricklayers' processes}, introduced in \cite{valak} and \cite{fluct} based on ideas of B.\ T\'oth. These models are slight generalizations of the misanthrope (see Cocozza-Thivent \cite{coco}) and also of the zero range processes. Each model from this class has a special one-parametered family of one dimensional discrete measures. If we fix that parameter and build the product of these measures for the sites, then we obtain a time-stationary distribution of the model. Hence we have a one-parametered family of equilibrium distributions; this parameter sets the \emph{slope of the wall}, a quantity corresponding to particle density in particle systems.

We consider product measures with the same marginals as for the stationary distribution, except for that we allow the parameter to change from site to site. The main result states that distributions in this class evolve to linear combinations of distributions from the same class, but only in the special case of \emph{exponential} bricklayers' process. The form of the linear coefficients allows us to interpret some of the situations as ordinary random walks of shocks having discontinuity of size one. We also obtain the nature of interaction between these shocks. It follows that a group formed by a number of such one-sized jump shocks is of stochastically bounded size in time, i.e.\ shocks of larger (integer-valued) discontinuities represented by such a group are sharp under any kind of hydrodynamic scaling.

Belitsky and Sch\"utz \cite{qse} derive results similar to ours for the simple exclusion process with the use of a quantum algebra symmetry. Unfortunately, we were not able to develop the quantum formalism for our locally infinite state space systems. However, there is a remarkable analogue between their work and the present settings: in \cite{qse} special relations are required between the shock densities and the particle jump rates. These relations are identical to those of Derrida, Lebowitz and Speer \cite{dls} allowing there an exact (not only asymptotic) product-description of the stationary distribution as seen from the second class particle. In the present work we shall be specially requiring integer values for the size of discontinuities of our shocks, in a very similar way as needed for the exact result of Bal\'azs \cite{valak} from the viewpoint of the second class particle.

The model and its stationary distributions are introduced in the next section. Then we give a very brief overview on the Eulerian hydrodynamic limit of these models in Section \ref{sc:hyd}. Our results are contained in Section \ref{sc:res}, and the proofs can be found in the last section. 

\fej{The bricklayers' model}\label{sc:model}

\alfej{Infinitesimal generator}

In \cite{fluct}, based on the idea of B.\ T\'oth, Bal\'azs defines a class of deposition models including the well-known simple exclusion and zero range processes. The so-called bricklayers' model to be discussed below is a family from the class described there. We briefly repeat the definition here, see \cite{fluct} for a more detailed introduction.

\begin{figure}[p]
\begin{center}
\begin{picture}(100, 160)(0, -10)
\linethickness{0.2pt}
\put(0, 10){\line(1, 0){100}}
\put(0, 30){\line(1, 0){100}}
\put(0, 50){\line(1, 0){100}}
\put(0, 70){\line(1, 0){30}}
\put(90, 70){\line(1, 0){10}}
\put(0, 90){\line(1, 0){10}}

\put(10, 0){\line(0, 1){90}}
\put(30, 0){\line(0, 1){70}}
\put(50, 0){\line(0, 1){50}}
\put(70, 0){\line(0, 1){50}}
\put(90, 0){\line(0, 1){50}}

\put(29, -7){$\scriptstyle{i}$}
\put(45, -7){$\scriptstyle{i\!+\!1}$}

\linethickness{2pt}

\put(0, 110){\line(1, 0){10}}
\put(10, 90){\line(1, 0){20}}
\put(30, 70){\line(1, 0){20}}
\put(50, 50){\line(1, 0){40}}
\put(90, 90){\line(1, 0){10}}

\put(10, 90){\line(0, 1){20}}
\put(30, 70){\line(0, 1){20}}
\put(50, 50){\line(0, 1){20}}
\put(90, 50){\line(0, 1){40}}

\put(32, 78){$\scriptstyle{\bigl\}\om_i}$}
\put(52, 58){$\scriptstyle{\bigl\}\om_{i+1}}$}

\end{picture}
\begin{picture}(110, 160)
\put(30, 60){\makebox(50, 12){$\underrightarrow{r(\om_i)+r(-\om_{i+1})\ }$}}
\end{picture}
\begin{picture}(100, 160)(0, -10)

\linethickness{0.2pt}

\put(0, 10){\line(1, 0){100}}
\put(0, 30){\line(1, 0){100}}
\put(0, 50){\line(1, 0){100}}
\put(0, 70){\line(1, 0){30}}
\put(90, 70){\line(1, 0){10}}
\put(0, 90){\line(1, 0){10}}

\put(10, 0){\line(0, 1){90}}
\put(30, 0){\line(0, 1){70}}
\put(50, 0){\line(0, 1){50}}
\put(70, 0){\line(0, 1){50}}
\put(90, 0){\line(0, 1){50}}

\put(29, -7){$\scriptstyle{i}$}
\put(45, -7){$\scriptstyle{i\!+\!1}$}

\linethickness{2pt}

\put(0, 110){\line(1, 0){10}}
\put(10, 90){\line(1, 0){40}}
\put(50, 50){\line(1, 0){40}}
\put(90, 90){\line(1, 0){10}}

\put(10, 90){\line(0, 1){20}}
\put(50, 50){\line(0, 1){40}}
\put(90, 50){\line(0, 1){40}}

\put(30, 70){\dashbox{2}(20, 20){}}

\end{picture}
\end{center}
\caption{A possible move}\label{fig:elso}
\end{figure}
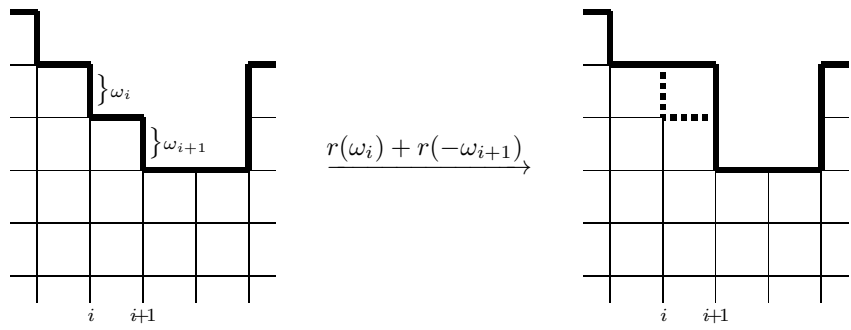

We consider the state space 
\begin{equation}
\Omega=\left\{\un\om=(\om_i)_{i\in\Zb}\ :\ \om_i\in\Zb\right\}=\Zb^{\Zb}.\label{eq:Omega}
\end{equation}
For each pair of neighboring sites $i$ and $i+1$ of $\Zb$, we can imagine a column built of bricks, above the edge $(i,\,i+1)$. The height of this column is denoted by $h_i$. If $\un{\om}(t)\in\Omega$ for a fixed time $t\in\Rb$ then $\om_i(t)=h_{i-1}(t)-h_i(t)\,\in\Zb$ is the negative discrete gradient of the height of the ``wall''. We are basically interested in this quantity, not the absolute heights $h$ of columns. Given a configuration, the growth of a column is described by a Poisson process. A brick can be added to a column, with rates depending on the local configuration:
\[
\un\om\longrightarrow\omi\ \ \text{with rate}\ \ r(\om_i)+r(-\om_{i+1}),
\]
where
\[
\omi=(\dots,\,\om_{i-1},\,\om_i-1,\,\om_{i+1}+1,\,\om_{i+2},\,\dots)
\]
only differs from $\un\om$ at sites $i$ and $i+1$. See fig.\,\ref{fig:elso} for some possible instantaneous changes. The process can be represented by bricklayers standing at each site $i$, laying a brick on the column on their right with rate $r(\om_i)$ and laying a brick to their left with rate $r(-\om_i)$. This interpretation gives reason to call these model bricklayers' model.
\afterpage{\clearpage}

The rate function $r$ is a positive real valued function on $\Zb$. In order to have product measures as stationary ones (see Corollary \ref{cor:Gibbs}), we assume that $r$ satisfies
\begin{equation}
r(z)\cdot r(-z+1)=1\qquad\text{for all}\ z\in\Zb.\label{eq:ratafelt}
\end{equation}

At time $t$, the interface mentioned before is described by $\un{\om}(t)$. Let $\varphi\,:\,\Omega\to\Rb$ be a bounded cylinder function i.e.\ $\varphi$ depends on a finite number of variables $\om_i$. The growth of this interface is a Markov process, with the formal infinitesimal generator $L$:
\begin{equation}
(L\varphi)(\un\om)=\sum_{i\in\Zb}\Bigl\{\left[r(\om_i)+r(-\om_{i+1})\right]\cdot\left[\vp(\omi)-\vp(\un\om)\right]\Bigr\}.\label{eq:gen}
\end{equation}
Note that for each index $i,\ \om_i$ can also be negative hence direct particle interpretation fails. Apart from this, $L$ is formally the sum of two totally asymmetric zero range-generators; one for $\om$ as number of particles which jump to the right, and another for $-\om$ as number of particles with jumps to the left. Therefore, part of our arguments can easily be applied to zero range. We shall notify the reader where the analogy ends. Unfortunately, the most interesting conclusions do not apply to the zero range process.

When constructing the process rigorously, problems may arise due to the unbounded growth rates. Existing methods for constructing such processes only apply when the rate function $r$ does not grow faster than linear, see Andjel \cite{and}, Booth and Quant \cite{lorna}. However, models with exponential rates $r$ are of essential interest for this note. In a forthcoming paper we prove that, in the attractive case, when $r$ is a monotone increasing function and does not grow faster than exponentially, the growth of any fixed column until any fixed moment is a random variable of finite second moment if we start the process from an appropriate set of initial configurations. Unfortunately we are not able to carry out a full construction of this Markov-process at this moment. We assume thus that existence of dynamics on a set of tempered configurations $\widetilde\Omega$ (i.e.\ configurations obeying some restrictive growth conditions) can be established. Technically we assume that $\widetilde\Omega$ is of full measure w.r.t.\ the product measures defined in Section \ref{sc:gibbs}. We do not deal with this question in the present paper.

\hop
{\bf The exponential bricklayers' model}

\bigskip

\noindent A special case of the models is the exponential bricklayers' model (EBL), where for $z\in\Zb$
\begin{equation}
r(z)=\e{-\frac{\be}{2}}\,\e{\be z}\label{eq:eblr}
\end{equation}
with a positive real parameter $\be$.

\alfej{The product measure}\label{sc:gibbs}

For $n\in\Zb^+$, we define 
\[
r(n)!:\,=\prod_{y=1}^nr(y)
\]
with the convention that the empty product has value one. For general $z\in\Zb$, we then define
\[
r(z)!:\,=r(|z|)!.
\]
Using \eqref{eq:ratafelt}, the property
\[
r(z+1)!=r(z)!\cdot r(z+1)
\]
can be shown for any integer $z$, giving justification of the notation used here. Let
\[
\bar\te:\,=\liminf_{n\to\infty}\log(r(n))
\]
which we assume to be strictly positive. By \eqref{eq:ratafelt}, it is in fact strictly positive in case of monotonicity of $r$ or, in other words, attractivity of the model. $\bar\te$ can even be infinite. With a generic real parameter $\te\in\left(-\bar\te,\,\bar\te\right)$, we define 
\[
Z(\te):\,=\sum_{z=-\infty}^\infty\frac{\e{\te z}}{r(z)!}
\]
and the measure
\begin{equation}
\mu^{(\te)}(z):\,=\frac{1}{Z(\te)}\cdot\frac{\e{\te z}}{r(z)!}\label{eq:om}
\end{equation}
on $\Zb$. The definition leads to two simple but essential properties:
\begin{equation}
\begin{split}
\e{\te}\cdot\frac{\mu^{(\te)}(z-1)}{\mu^{(\te)}(z)}&=r(z)\\
\e{-\te}\cdot\frac{\mu^{(\te)}(z+1)}{\mu^{(\te)}(z)}&=r(-z),
\end{split}\label{eq:miap}
\end{equation}
of which the latter is a consequence of \eqref{eq:ratafelt}. Based on these, we also see
\begin{equation}
\Ev\,r(z)=\e{\te}\quad\text{and}\quad\Ev\,r(-z)=\e{-\te}\label{eq:er}
\end{equation}
w.r.t.\ $\mu^{(\te)}$.

For the special case of the EBL model, for $\te\in(-\infty,\,\infty)$, we obtain the discrete normal distribution
\begin{equation}
\mu^{(\te)}(z)=\frac{\e{-\frac{\be}{2}\left(z-\frac{\te}{\be}\right)^2}}{\e{-\frac{\te^2}{2\be}}Z(\te)}=\frac{\e{-\frac{\be}{2}(z-m)^2}}{\widetilde{Z}(\be,\,m)}\label{eq:emu}
\end{equation}
with the notation $m:=\te/\be$.

For the purposes of this paper, we introduce the parameter-vector
\[
\un\te:\,=(\te_i)_{i\in\Zb}\quad\text{where}\ \te_i\in(-\bar\te,\,\bar\te),
\]
and the product distribution $\un\mu^{(\un\te)}$ on $\Omega$, having marginal $\mu^{(\te_i)}$ for site $i$. We also assume that $\te_i$ has finitely many different values as $i$ goes from $-\infty$ to $\infty$.

\fej{Hydrodynamic limit}\label{sc:hyd}

As we shall see later, the product measure $\un\mu^{(\un\te)}$ is stationary for the models if $\te_i$ has the same value for all sites $i$. Allowing this equilibrium to change on a macroscopic level, the hydrodynamic limit of these models can be derived. For different types of scaling limits we refer to T\'oth and Valk\'o \cite{hydro}. Here we only consider the so-called Euler-scaling. Strictly speaking, as in case of construction of these processes, we only know about the case of sublinear $r$ functions handled in the context of hydrodynamic limits. However, for this note exponential growth rates are essential. Via formal computations, one can guess for the form of the Eulerian limiting hydrodynamic partial differential equation even in this superlinear case. The variable of this equation is formally defined by
\begin{equation}
u(t,\,x):\,=\Ev\,\om_{x/\ve}(t/\ve).\label{eq:miau}
\end{equation}
Then via formal computations we obtain the differential equation
\begin{equation}
\frac{\pt u}{\pt t}+\frac{\pt J(u)}{\pt x}=0\label{eq:pde}
\end{equation}
as $\ve\to0$, where $J(u)$ is defined as follows. The function $u(\te)=\Ev^{(\te)}(\om)$ of $\te$ is strictly increasing since the derivative
\[
\frac{\di u(\te)}{\di\te}=\left(\Ev^{(\te)}(\om^2)-\left(\Ev^{(\te)}(\om)\right)^2\right)
\]
is positive ($-\bar\te<\te<\bar\te$). Let $\te(u)$ be the inverse function. The quantity $\Ev^{(\te)}\left(r(\om)+r(-\om)\right)$ depends on $\te$, and $J$ is defined by
\begin{equation}
J(u):=\Ev^{(\te(u))}\left(r(\om)+r(-\om)\right)=2\cosh(\te(u)).\label{eq:miaj}
\end{equation}

In \cite{fluct}, Bal\'azs proves that this function is strictly convex if $r$ is monotone increasing and is convex. (An example of such rates is the exponential bricklayers' process defined above.) In this case, \eqref{eq:pde} develops shock (weak) solutions. A (single) shock solution starts with initial data 
\[
u(0,\,x)=\left\{\begin{array}{lcl}
\!\!\!u_{\text{left}}&,\ \ &x<0\\
\!\!\!u_{\text{right}}&,\ \ &x\ge0
\end{array}\right. 
\]
with $u_{\text{left}}>u_{\text{right}}$. The stable weak solution is of the form
\[
u(t,\,x)=\left\{\begin{array}{lcl}
\!\!\!u_{\text{left}}&,\ \ &x<st\\
\!\!\!u_{\text{right}}&,\ \ &x\ge st
\end{array}\right. 
\]
where the speed $s$ of the traveling shock is determined by the Rankine-Hugoniot formula
\begin{equation}
s=\frac{J(u_{\text{left}})-J(u_{\text{right}})}{u_{\text{left}}-u_{\text{right}}}\ \ ,\label{eq:ran}
\end{equation}
see e.g.\ Smoller \cite{smo}. We define the \emph{size} of the shock as $u_\text{left}-u_\text{right}$.

In the case of multiple (but finitely many) shocks, i.e.\ piecewise constant, decreasing initial data with more than one discontinuities, a similar formula holds for each shock if $u$ has values $u_\text{left}$ and $u_\text{right}$ at its two sides. By convexity of $J$, these shocks meet in finite time, forming less and less shocks of larger and larger sizes.

Using definitions \eqref{eq:miaj}, \eqref{eq:miau} and the fact \eqref{eq:er} in the Rankine-Hugoniot formula \eqref{eq:ran}, the speed of a shock can now be written as
\begin{equation}
s=\frac{\e{\te(u_{\text{left}})}+\e{-\te(u_{\text{left}})}-\e{\te(u_{\text{right}})}-\e{-\te(u_{\text{right}})}}{u_{\text{left}}-u_{\text{right}}}.\label{eq:mitran}
\end{equation}

This is what we see on the macroscopic scale. Our aim now is to find some similar structures in the microscopic level.

\fej{Results}\label{sc:res}

For simplicity, we introduce the notations
\begin{equation}
\begin{array}{r@{\,}l}
\teip:&=(\dots,\,\te_{i-1},\,\te_i+\be,\,\te_{i+1},\,\dots)\\
\teim:&=(\dots,\,\te_{i-1},\,\te_i-\be,\,\te_{i+1},\,\dots)
\end{array}\label{eq:ipm}
\end{equation}
given a parameter $\be$.

\begin{tm}\label{tm:fo}
For the exponential and no other bricklayers' model, (neither for any totally asymmetric nearest-neighbor zero range process), the product measure $\un\mu^{(\un\te)}$, defined at the very end of Section \ref{sc:model}, evolves to a linear combination of product measures with similar structure. The evolution is formally described by the equation
\begin{multline}
\frac{\di}{\di t}\un\mu^{(\un\te)}=\sum_{i:\,\te_i\ne\te_{i+1}}\left(\e{\te_i}-\e{\te_{i+1}}\right)\cdot\left(\un\mu^{(\teipp)}-\un\mu^{(\un\te)}\right)+\\
+\sum_{i:\,\te_{i-1}\ne\te_i}\left(\e{-\te_i}-\e{-\te_{i-1}}\right)\cdot\left(\un\mu^{(\teimm)}-\un\mu^{(\un\te)}\right).\label{eq:muder}
\end{multline}
\end{tm}

The proof can be found in Section \ref{sc:pr}. We cannot give a nice interpretation of this evolution in the general setting. However, in some special cases, this theorem gives rise to a rather interesting observation. For the rest of this section, we are talking about the case of exponential bricklayers' process.

\alfej{One single (and special) shock}\label{sc:osing}

Fix a parameter $\te_\text{left}$ for the left-hand side, and $\te_\text{right}:\,=\te_\text{left}-\be$. Fix a position $Q\in\Zb$ and
\[
\te_i:\,=\left\{\begin{array}{ll}\te_\text{left}&\text{for}\ i\le Q-1,\\\te_\text{right}&\text{for}\ i\ge Q\end{array}\right.
\]
for components of $\un\te$. Then $\un\mu^{(\un\te)}$ is a shock-measure corresponding to a shock of size one: for the exponential model, a change of size $\be$ in the parameter $\te$ implies a jump of size one in the expectation $u=\Ev^{(\te)}(\om)$ (see also Lemma \ref{lm:eltol}). We say that this shock-measure \emph{is located at} position $Q-\frac12$. Bal\'azs \cite{valak} shows that this measure is in fact stationary as seen from the random position of the so-called \emph{second class particle}. In the present setting, we examine this structure from a fixed position. First, we apply \eqref{eq:muder} to this $\un\mu^{(\un\te)}$:
\begin{multline}
\frac{\di}{\di t}\un\mu^{(\un\te)}=\left(\e{\te_\text{left}}-\e{\te_\text{right}}\right)\cdot\left(\un\mu^{(\un\te^{(Q,\,+)})}-\un\mu^{(\un\te)}\right)+\\
+\left(\e{-\te_\text{right}}-\e{-\te_\text{left}}\right)\cdot\left(\un\mu^{(\un\te^{(Q-1,\,-)})}-\un\mu^{(\un\te)}\right).\label{eq:egysok}
\end{multline}
Then we realize that
\[
\te^{(Q,\,+)}_i=\left\{\begin{array}{ll}\te_\text{left}&\text{for}\ i\le Q,\\\te_\text{right}&\text{for}\ i\ge Q+1,\end{array}\right.
\]
and
\[
\te^{(Q-1,\,-)}_i=\left\{\begin{array}{ll}\te_\text{left}&\text{for}\ i\le Q-2,\\\te_\text{right}&\text{for}\ i\ge Q-1.\end{array}\right.
\]
The former is the same shock measure located at $Q+\frac12$, the latter is located at $Q-1\frac12$, see figures \ref{fig:masodik} and \ref{fig:harmadik}.
\begin{figure}[p]
\begin{center}
\begin{picture}(100, 120)(-10,-40)
\linethickness{0.2pt}
\put(0,10){\vector(1,0){100}}
\put(10,0){\vector(0,1){70}}
\multiput(30,7)(20,0){4}{\line(0,1){6}}
\multiput(7,20)(0,20){3}{\line(1,0){6}}
\multiput(30,15)(0,5){9}{\line(0,1){2}}
\multiput(50,15)(0,5){9}{\line(0,1){2}}
\multiput(70,15)(0,5){5}{\line(0,1){2}}
\multiput(90,15)(0,5){5}{\line(0,1){2}}
\put(58,48){$\Bigl\}\be$}
\put(100,0){$i$}
\put(0,68){$\te_i$}
\put(41,0){$\scriptstyle{Q-1}$}
\put(67,0){$\scriptstyle{Q}$}
\put(-10,57){$\scriptstyle{\te_\text{left}}$}
\put(-10,37){$\scriptstyle{\te_\text{right}}$}

\put(0,-25){\vector(1,0){100}}
\multiput(20,-27)(20,0){4}{\line(0,1){4}}
\put(57,-20){$\bullet$}
\put(97,-33){$\scriptstyle{\Zb+\frac12}$}
\put(20,-35){\footnotesize{the random walker}}

\linethickness{1pt}
\put(22,60){\line(1,0){16}}
\put(42,60){\line(1,0){16}}
\put(62,40){\line(1,0){16}}
\put(82,40){\line(1,0){16}}

\end{picture}
\begin{picture}(100, 120)
\put(25, 75){\makebox(50, 12){$\underrightarrow{\e{\te_\text{left}}-\e{\te_\text{right}}}$}}
\end{picture}
\begin{picture}(100, 120)(0,-40)
\linethickness{0.2pt}
\put(0,10){\vector(1,0){100}}
\put(10,0){\vector(0,1){70}}
\multiput(30,7)(20,0){4}{\line(0,1){6}}
\multiput(7,20)(0,20){3}{\line(1,0){6}}
\multiput(30,15)(0,5){9}{\line(0,1){2}}
\multiput(50,15)(0,5){9}{\line(0,1){2}}
\multiput(70,15)(0,5){9}{\line(0,1){2}}
\multiput(90,15)(0,5){5}{\line(0,1){2}}
\put(78,48){$\Bigl\}\be$}
\put(100,0){$i$}
\put(0,68){$\te_i$}
\put(41,0){$\scriptstyle{Q-1}$}
\put(67,0){$\scriptstyle{Q}$}
\put(-10,57){$\scriptstyle{\te_\text{left}}$}
\put(-10,37){$\scriptstyle{\te_\text{right}}$}

\put(0,-25){\vector(1,0){100}}
\multiput(20,-27)(20,0){4}{\line(0,1){4}}
\put(77,-20){$\bullet$}
\put(97,-33){$\scriptstyle{\Zb+\frac12}$}
\put(20,-35){\footnotesize{the random walker}}
\put(62, -20){$\dashrightarrow$}

\linethickness{1pt}
\put(22,60){\line(1,0){16}}
\put(42,60){\line(1,0){16}}
\put(62,60){\line(1,0){16}}
\put(82,40){\line(1,0){16}}

\end{picture}
\end{center}
\caption{Step of the random walker from $Q-\frac12$ to $Q+\frac12$}\label{fig:masodik}
\end{figure}%
\afterpage{\clearpage}%
\begin{figure}[p]
\begin{center}
\begin{picture}(100, 120)(-10,-40)
\linethickness{0.2pt}
\put(0,10){\vector(1,0){100}}
\put(10,0){\vector(0,1){70}}
\multiput(30,7)(20,0){4}{\line(0,1){6}}
\multiput(7,20)(0,20){3}{\line(1,0){6}}
\multiput(30,15)(0,5){9}{\line(0,1){2}}
\multiput(50,15)(0,5){9}{\line(0,1){2}}
\multiput(70,15)(0,5){5}{\line(0,1){2}}
\multiput(90,15)(0,5){5}{\line(0,1){2}}
\put(58,48){$\Bigl\}\be$}
\put(100,0){$i$}
\put(0,68){$\te_i$}
\put(41,0){$\scriptstyle{Q-1}$}
\put(67,0){$\scriptstyle{Q}$}
\put(-10,57){$\scriptstyle{\te_\text{left}}$}
\put(-10,37){$\scriptstyle{\te_\text{right}}$}

\put(0,-25){\vector(1,0){100}}
\multiput(20,-27)(20,0){4}{\line(0,1){4}}
\put(57,-20){$\bullet$}
\put(97,-33){$\scriptstyle{\Zb+\frac12}$}
\put(20,-35){\footnotesize{the random walker}}

\linethickness{1pt}
\put(22,60){\line(1,0){16}}
\put(42,60){\line(1,0){16}}
\put(62,40){\line(1,0){16}}
\put(82,40){\line(1,0){16}}

\end{picture}
\begin{picture}(100, 120)
\put(25, 75){\makebox(50, 12){$\underrightarrow{\e{-\te_\text{right}}-\e{-\te_\text{left}}}$}}
\end{picture}
\begin{picture}(100, 120)(0,-40)
\linethickness{0.2pt}
\put(0,10){\vector(1,0){100}}
\put(10,0){\vector(0,1){70}}
\multiput(30,7)(20,0){4}{\line(0,1){6}}
\multiput(7,20)(0,20){3}{\line(1,0){6}}
\multiput(30,15)(0,5){9}{\line(0,1){2}}
\multiput(50,15)(0,5){5}{\line(0,1){2}}
\multiput(70,15)(0,5){5}{\line(0,1){2}}
\multiput(90,15)(0,5){5}{\line(0,1){2}}
\put(38,48){$\Bigl\}\be$}
\put(100,0){$i$}
\put(0,68){$\te_i$}
\put(41,0){$\scriptstyle{Q-1}$}
\put(67,0){$\scriptstyle{Q}$}
\put(-10,57){$\scriptstyle{\te_\text{left}}$}
\put(-10,37){$\scriptstyle{\te_\text{right}}$}

\put(0,-25){\vector(1,0){100}}
\multiput(20,-27)(20,0){4}{\line(0,1){4}}
\put(37,-20){$\bullet$}
\put(97,-33){$\scriptstyle{\Zb+\frac12}$}
\put(20,-35){\footnotesize{the random walker}}
\put(42, -20){$\dashleftarrow$}

\linethickness{1pt}
\put(22,60){\line(1,0){16}}
\put(42,40){\line(1,0){16}}
\put(62,40){\line(1,0){16}}
\put(82,40){\line(1,0){16}}

\end{picture}
\end{center}
\caption{Step of the random walker from $Q-\frac12$ to $Q-1\frac12$}\label{fig:harmadik}
\end{figure}%
\afterpage{\clearpage}%
Thus \eqref{eq:egysok} says the following: a shock measure of size one, located at $Q-\frac12$ becomes a similar measure at $Q+\frac12$ with rate $\e{\te_\text{left}}-\e{\te_\text{right}}$, or becomes one at $Q-1\frac12$ with rate $\e{-\te_\text{right}}-\e{-\te_\text{left}}$. Hence the shock performs an ordinary nearest-neighbor random walk with constant rates $\e{\te_\text{left}}-\e{\te_\text{right}}$ to the left and $\e{-\te_\text{right}}-\e{-\te_\text{left}}$ to the right, respectively. These rates are positive according to $\te_\text{left}>\te_\text{right}$, and they in fact agree with the \emph{expectation} of the left and right jump rates of the second class particle w.r.t.\ $\un\mu^{(\un\te)}$. The expected speed of this random walk therefore agrees with the speed of the macroscopic shock (of size one) suggested by the Rankine-Hugoniot formula \eqref{eq:mitran}. This is the shock-measure of \cite{valak} as seen from outside.

We emphasize that the ordinary random walk shown here and the second class particle of \cite{valak} are different objects. While the jump rates of the latter depend on the actual state of the model, in the present context we only deal with distributions of the process, hence we can not even state any kind of dependence of our random walks on particular states of the model. One could possibly think about comparing the law of the second class particle, integrated out w.r.t.\ the initial distribution $\un\mu^{(\un\te)}$ of the model, to the distribution of our random walkers. We do not know about any argument which could bring these laws in connection to each other. In particular, by averaging out the continuously changing random environment this way, nothing assures {\sl a priori} the Markov-property of the second class particle, while our ordinary random walkers clearly have this behavior.

\alfej{Multiple (special) shocks}\label{sc:muso}

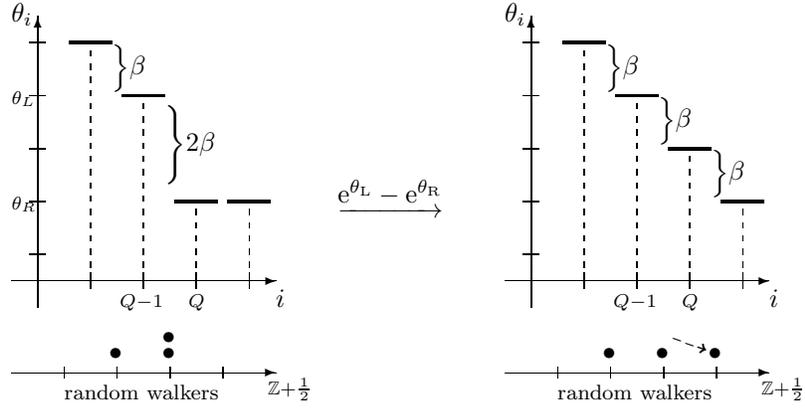
\begin{figure}[p]
\begin{center}
\begin{picture}(100, 150)(0,-40)
\linethickness{0.2pt}
\put(0,10){\vector(1,0){100}}
\put(10,0){\vector(0,1){110}}
\multiput(30,7)(20,0){4}{\line(0,1){6}}
\multiput(7,20)(0,20){5}{\line(1,0){6}}
\multiput(30,15)(0,5){17}{\line(0,1){2}}
\multiput(50,15)(0,5){13}{\line(0,1){2}}
\multiput(70,15)(0,5){5}{\line(0,1){2}}
\multiput(90,15)(0,5){5}{\line(0,1){2}}
\put(38,88){$\Bigl\}\be$}
\put(58,59){$\Biggl\}2\be$}
\put(100,0){$i$}
\put(0,108){$\te_i$}
\put(41,0){$\scriptstyle{Q-1}$}
\put(67,0){$\scriptstyle{Q}$}
\put(0,77){$\scriptstyle{\te_L}$}
\put(0,37){$\scriptstyle{\te_R}$}

\put(0,-25){\vector(1,0){100}}
\multiput(20,-27)(20,0){4}{\line(0,1){4}}
\multiput(37,-20)(20,0){2}{$\bullet$}
\put(57,-14){$\bullet$}
\put(97,-33){$\scriptstyle{\Zb+\frac12}$}
\put(20,-35){\footnotesize{random walkers}}

\linethickness{1pt}
\put(22,100){\line(1,0){16}}
\put(42,80){\line(1,0){16}}
\put(62,40){\line(1,0){16}}
\put(82,40){\line(1,0){16}}

\end{picture}
\begin{picture}(80, 150)
\put(15, 75){\makebox(50, 12){$\underrightarrow{\e{\te_\text{L}}-\e{\te_\text{R}}}$}}
\end{picture}
\begin{picture}(100, 150)(0,-40)
\linethickness{0.2pt}
\put(0,10){\vector(1,0){100}}
\put(10,0){\vector(0,1){110}}
\multiput(30,7)(20,0){4}{\line(0,1){6}}
\multiput(7,20)(0,20){5}{\line(1,0){6}}
\multiput(30,15)(0,5){17}{\line(0,1){2}}
\multiput(50,15)(0,5){13}{\line(0,1){2}}
\multiput(70,15)(0,5){9}{\line(0,1){2}}
\multiput(90,15)(0,5){5}{\line(0,1){2}}
\multiput(38,88)(20,-20){3}{$\Bigl\}\be$}
\put(100,0){$i$}
\put(0,108){$\te_i$}
\put(41,0){$\scriptstyle{Q-1}$}
\put(67,0){$\scriptstyle{Q}$}

\put(0,-25){\vector(1,0){100}}
\multiput(20,-27)(20,0){4}{\line(0,1){4}}
\multiput(37,-20)(20,0){3}{$\bullet$}
\put(97,-33){$\scriptstyle{\Zb+\frac12}$}
\put(20,-35){\footnotesize{random walkers}}
\put(62, -14){\rotatebox{-19}{$\dashrightarrow$}}

\linethickness{1pt}
\put(22,100){\line(1,0){16}}
\put(42,80){\line(1,0){16}}
\put(62,60){\line(1,0){16}}
\put(82,40){\line(1,0){16}}

\end{picture}
\end{center}
\caption{A shock of size one steps to the right from position $Q-\frac12$}\label{fig:negyedik}
\end{figure}

\begin{figure}[p]
\begin{center}
\begin{picture}(100, 150)(0,-40)
\linethickness{0.2pt}
\put(0,10){\vector(1,0){100}}
\put(10,0){\vector(0,1){110}}
\multiput(30,7)(20,0){4}{\line(0,1){6}}
\multiput(7,20)(0,20){5}{\line(1,0){6}}
\multiput(30,15)(0,5){17}{\line(0,1){2}}
\multiput(50,15)(0,5){13}{\line(0,1){2}}
\multiput(70,15)(0,5){5}{\line(0,1){2}}
\multiput(90,15)(0,5){5}{\line(0,1){2}}
\put(38,88){$\Bigl\}\be$}
\put(58,59){$\Biggl\}2\be$}
\put(100,0){$i$}
\put(0,108){$\te_i$}
\put(41,0){$\scriptstyle{Q-1}$}
\put(67,0){$\scriptstyle{Q}$}
\put(0,77){$\scriptstyle{\te_L}$}
\put(0,37){$\scriptstyle{\te_R}$}

\put(0,-25){\vector(1,0){100}}
\multiput(20,-27)(20,0){4}{\line(0,1){4}}
\multiput(37,-20)(20,0){2}{$\bullet$}
\put(57,-14){$\bullet$}
\put(97,-33){$\scriptstyle{\Zb+\frac12}$}
\put(20,-35){\footnotesize{random walkers}}

\linethickness{1pt}
\put(22,100){\line(1,0){16}}
\put(42,80){\line(1,0){16}}
\put(62,40){\line(1,0){16}}
\put(82,40){\line(1,0){16}}

\end{picture}
\begin{picture}(80, 150)
\put(15, 75){\makebox(50, 12){$\underrightarrow{\e{-\te_\text{R}}-\e{-\te_\text{L}}}$}}
\end{picture}
\begin{picture}(100, 150)(0,-40)
\linethickness{0.2pt}
\put(0,10){\vector(1,0){100}}
\put(10,0){\vector(0,1){110}}
\multiput(30,7)(20,0){4}{\line(0,1){6}}
\multiput(7,20)(0,20){5}{\line(1,0){6}}
\multiput(30,15)(0,5){17}{\line(0,1){2}}
\multiput(50,15)(0,5){9}{\line(0,1){2}}
\multiput(70,15)(0,5){5}{\line(0,1){2}}
\multiput(90,15)(0,5){5}{\line(0,1){2}}
\put(58,48){$\Bigl\}\be$}
\put(38,79){$\Biggl\}2\be$}
\put(100,0){$i$}
\put(0,108){$\te_i$}
\put(41,0){$\scriptstyle{Q-1}$}
\put(67,0){$\scriptstyle{Q}$}

\put(0,-25){\vector(1,0){100}}
\multiput(20,-27)(20,0){4}{\line(0,1){4}}
\multiput(37,-20)(20,0){2}{$\bullet$}
\put(37,-14){$\bullet$}
\put(97,-33){$\scriptstyle{\Zb+\frac12}$}
\put(20,-35){\footnotesize{random walkers}}
\put(42, -14){$\dashleftarrow$}

\linethickness{1pt}
\put(22,100){\line(1,0){16}}
\put(42,60){\line(1,0){16}}
\put(62,40){\line(1,0){16}}
\put(82,40){\line(1,0){16}}

\end{picture}
\end{center}
\caption{A shock of size one steps to the left from position $Q-\frac12$}\label{fig:otodik}
\end{figure}
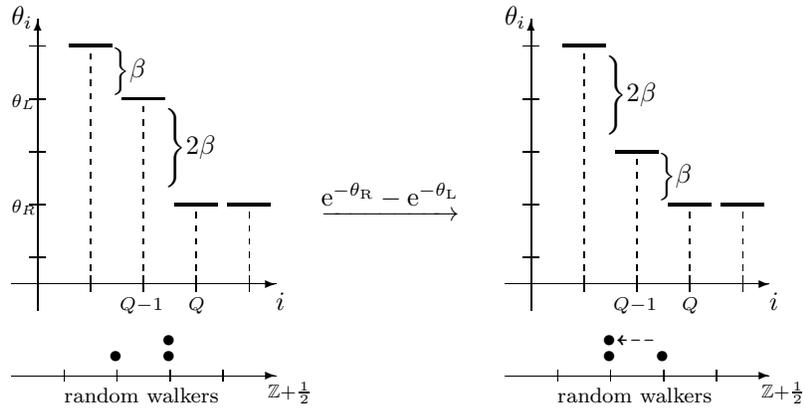

As for coordinates of $\un\te$, assume that $\te_i$ as a function of $i$
\begin{itemize}
\item is decreasing,
\item changes finitely many times
\item in multiples of $\be$.
\end{itemize}
We also assume that the values of $\te_i$ range from $\te_\text{right}$ to $\te_\text{left}=\te_\text{right}+n\be$. As each jump of size $\be$ in the parameter $\te$ results in a jump of size one in the expectation (see Lemma \ref{lm:eltol}), and so in the hydrodynamic quantity $u$, $\un\mu^{(\un\te)}$ now corresponds to a multiple shock measure, where each shock has integer size giving a sum of $n$. In a similar way as for the single-shock case, one reads the following representation of \eqref{eq:muder}. If $\te_\text{R}=\te_\text{L}-m\be$, and
\begin{equation}
\te_i=\left\{\begin{array}{ll}\te_\text{L}&\text{for}\ i=Q-1,\\\te_\text{R}&\text{for}\ i=Q,\end{array}\right.\label{eq:LR}
\end{equation}
then we have a shock of size $m$ located at $Q-\frac12$. With rate $\e{\te_\text{L}}-\e{\te_\text{R}}$, $\te_Q$ increases by $\be$, i.e.\ the size of the shock at $Q-\frac12$ decreases by one and the size of the shock at $Q+\frac12$ increases by one (or a new shock of size one emerges if there wasn't any shock here before). With rate $\e{-\te_\text{R}}-\e{-\te_\text{L}}$, $\te_{Q-1}$ decreases by $\be$, i.e.\ the size of the shock at $Q-\frac12$ decreases by one and the size of the shock at $Q-1\frac12$ increases by one (or a new shock of size one emerges here). See figures \ref{fig:negyedik} and \ref{fig:otodik} for these steps.
\afterpage{\clearpage}%
\afterpage{\clearpage}%
This, in fact, gives us the following picture: the $n$ one-sized shocks can be represented by $n$ random walkers, each performing a nearest-neighbor random walk on $\Zb+\frac12$. These walkers meet eventually, hence any number $m\le n$ of them can be present at a given position. Assume we have $m>0$ walkers at this given position, and $k$ of these walkers strictly left to this position. In this case, since each walker (i.e.\ each shock) changes the value of $\te$ by $\be$, we have $\te_\text{L}=\te_\text{left}-k\be$ and $\te_\text{R}=\te_\text{left}-(k+m)\be$ with the notation of \eqref{eq:LR}. Therefore, the rate with which any walker steps from here is%
{\renewcommand\arraystretch{1.4}
\begin{equation}
\begin{array}{r@{\,}c@{}cl}
\e{\te_\text{L}}-\e{\te_\text{R}}=&\e{\te_\text{left}-k\be}&\cdot\left(1-\e{-m\be}\right)&\text{to the right and}\\
\e{-\te_\text{R}}-\e{-\te_\text{L}}=&\e{-\te_\text{left}+k\be}&\cdot\left(\e{m\be}-1\right)&\text{to the left,}
\end{array}\label{eq:sokLR}
\end{equation}}%
respectively. If a walker is alone at her site and she has $l\ge0$ walkers to her left, then we have $m=1$ and write $l$ instead of $k$. In this case, our formulas reduce to $\e{\te_\text{left}-l\be}\left(1-\e{-\be}\right)$ as for rate to jump to the right with, and $\e{-\te_\text{left}+l\be}\left(\e{\be}-1\right)$ as for rate to jump to the left. In fact, the jump rates \eqref{eq:sokLR} are just the sums of the latter jump rates for the single walkers as $l$ goes from $k$ to $m+k-1$. In other words, the rate with which any walker jumps from a position is just the sum of those rates if the walkers at that position were alone.

The expressions for these rates show that walkers on the left (smaller $k$ values) have larger rates to jump to the right and smaller rates to jump to the left than others more on the right (larger $k$ values). Only their order matters, not their absolute position. \emph{The form of these jump rates is the only way the shocks interact.}

We also see from here that the shock size at any fixed position can only change by one at a time. Hence shocks of size larger than one can never move, only shocks of size one are separated and merged again eventually somewhere else. \emph{On the microscopic level, larger shocks are only formed by accidental meeting of shocks of size one which perform the random walks}.

The rate of these walkers (or shocks) agree with the \emph{expected} jump rates of the second class particle w.r.t.\ $\un\mu^{(\un\te)}$. However, as mentioned above for the case of a single shock, the second class particle and our random walkers are conceptually different objects.

As larger shocks can only move by separation and displacement of one-sized shocks, one can ask if these several little shocks in fact show the properties of one large shock on the macroscopic scale. An answer to this question is given below: the following lemma shows that these type of shocks are sharp under any kind of hydrodynamic limit. We define the \emph{width} of such a larger shock (formed by a number of one-sized shocks) as the difference between the position of the leftmost and the rightmost shock (or random walker).
\begin{lm}\label{lm:w} 
The width of our types of shocks has bounded moments in time.
\end{lm}
\noindent
The proof of this lemma is contained in Section \ref{sc:pr}.

As the one-sized shocks (i.e.\ walkers) stay in stochastically bounded distance from each other, a natural quantity indicating the position of this group of walkers is their central mass. To be more precise, let $X^{(m)}(t)$ denote the position of the $m$-th walker at time $t$, $m=1\dots n$ (with any labeling), and define
\begin{equation}
X(t):\,=\frac{X^{(1)}(t)+X^{(2)}(t)+\dots+X^{(n)}(t)}{n}.\label{eq:cm}
\end{equation}
\begin{lm}\label{lm:cm}
The process $X(t)$ is an ordinary continuous-time random walk on $(\Zb+1/2)/n$ with
\[
\begin{array}{rl}
\text{right jump rate}&\e{\te_\text{left}}-\e{\te_\text{right}},\\
\text{left jump rate}&\e{-\te_\text{right}}-\e{-\te_\text{left}}.
\end{array}
\]
\end{lm}
\noindent
The proof of this lemma can also be found in Section \ref{sc:pr}. As a corollary, $X(t)$ satisfies strong law of large numbers
\[
\lim_{t\to\infty}\frac{X(t)}{t}=\frac{\e{\te_\text{left}}+\e{-\te_\text{left}}-\e{\te_\text{right}}-\e{-\te_\text{right}}}{n}\qquad\text{a.s.}
\]
As the $n$ walkers together represent a shock of jump size $n$, this formula is in accordance with the Rankine-Hugoniot formula \eqref{eq:mitran}, derived for the speed of shock solutions of the PDE.

The model and the situations detailed above are the ones for which we could establish existence of these kind of solutions in terms of ordinary random walks. Situations when the random walk interpretation fails are listed below.

\alfej{Other initial conditions}

Unfortunately, we don't have such a nice interpretation when any of the conditions in Section \ref{sc:muso} fails. Theorem \ref{tm:fo} still holds in this case as it does not depend on these conditions, hence the product structure of the measure is still closed under the time-evolution of the system. One just cannot give the random walk representation of the result. 
If $\te_i$ as function of $i$ is not decreasing at some $i$ value, then the rates of our random walkers would become negative; we don't know how to interpret this case. It is not surprising that such random walk description fails: upward jumps of the function $\te_i$ correspond to upward shock solutions. By strict convexity of the hydrodynamic flux function $J(u)$, which is in fact established in \cite{fluct} for the exponential rate functions $r$ we have in the exponential model, upward shocks are non stable, they evolve into rarefaction waves (see e.g.\ \cite{smo}). Hence it is not reasonable to expect some kind of random walks of upward jump solutions.

On the other hand, if $\te_i$ is decreasing but not in multiples of $\be$, that means that we have shocks of non-integer size (see Lemma \ref{lm:eltol}). By separation of one-sized shocks, as described by the Theorem, eventually a shock of size smaller than one appears. When another shock of size one jumps from its position, it leaves an upward shock behind itself. At this moment, our walkers' rates become negative, and the interpretation fails again. Unfortunately we don't have an intuitive argument explaining this phenomenon.

\alfej{The not completely asymmetric case}

Bricklayers' process in the presented form corresponds to totally asymmetric particle systems. One can modify the generator by multiplying the rates so far introduced with a factor $p$ and include brick-removals as well:
\[
\un\om\longrightarrow\omir\ \ \text{with rate}\ \ q\left(r(-\om_i)+r(\om_{i+1})\right),
\]
where
\[
\omir=(\dots,\,\om_{i-1},\,\om_i+1,\,\om_{i+1}-1,\,\om_{i+2},\,\dots)
\]
and $q+p=1$. Theorem \ref{tm:fo} can easily be extended to this case: the class of product measures with marginals $\mu^{(\te)}$ of spatially changing $\te$ parameters is closed under the time-evolution of the process. The question is if the random walk interpretation can also be applied to this more general case. While well-posed asymmetry of simple exclusion is essential in \cite{qse}, it proves to be useless for our random walk interpretation in the bricklayers' process. Independently of the initial distribution (either increasing, or decreasing), general asymmetry of the model would immediately imply negative jump rates of our random walkers when trying to introduce them. Hence we cannot generalize our random walk interpretation, or weaken any of the conditions of Section \ref{sc:muso}, needed for this interpretation, by introducing general asymmetry in the system.

\alfej{Other models}

As seen above, other bricklayers' models do not even show closedness of the product measures under the time-evolution as exponential bricklayers' process does. For general bricklayers' models or shock data, we do not expect qualitatively different behavior, rather results similar to the exact ones obtained here to hold asymptotically in space. For example, based on results obtained in this direction for simple exclusion (see e.g. Derrida, Lebowitz and Speer \cite{dls} or Ferrari, Fontes and Kohayakawa \cite{imes}), we conjecture that the stationary distribution as seen from the single one second class particle converges asymptotically to the product measure $\un\mu^{(\un\te)}$ with the appropriate $\te_\text{left}$ and $\te_\text{right}$ values. We do not have any results in this direction. There is only an easy argument in Section 4.2 of Bal\'azs \cite{valak} showing that once this convergence and also ergodicity of the stationary measure is established, the second class particle satisfies the Law of Large Numbers with the speed value given by the Rankine-Hugoniot condition.

\fej{Proofs}\label{sc:pr}

For simplicity, we introduce the notations
\begin{equation}
\begin{array}{r@{\,}l}
\omip:&=(\dots,\,\om_{i-1},\,\om_i+1,\,\om_{i+1},\,\dots)\\
\omim:&=(\dots,\,\om_{i-1},\,\om_i-1,\,\om_{i+1},\,\dots).
\end{array}\label{eq:ipmo}
\end{equation}
\begin{lm}\label{lm:fo}
Let $\vp$ be a bounded cylinder function on $\Omega$. Then according to the bricklayers' dynamics with any rate function $r$ (satisfying \eqref{eq:ratafelt}), 
\begin{multline}
\frac{\di}{\di t}\Evt\,\vp(\un\om)=\sum_{i:\,\te_i\ne\te_{i+1}}\left(\e{\te_i}-\e{\te_{i+1}}\right)\cdot\left(\Evt\,\vp(\omipp)-\Evt\,\vp(\un\om)\right)+\\
+\sum_{i:\,\te_{i-1}\ne\te_i}\left(\e{-\te_i}-\e{-\te_{i-1}}\right)\cdot\left(\Evt\,\vp(\omimm)-\Evt\,\vp(\un\om)\right),\label{eq:fo}
\end{multline}
where the expectation $\Evt$ is taken w.r.t.\ the product measure $\un\mu^{(\un\te)}$ defined at the end of Section \ref{sc:model}.
\end{lm}
\begin{proof}
Using the generator \eqref{eq:gen}, we start computing this time derivative. By \eqref{eq:miap} we get to
\begin{multline*}
\frac{\di}{\di t}\Evt\,\vp(\un\om)=\Evt\,(L\vp)(\un\om)=\\
=\Evt\left\{\sum_{i\in\Zb}\left(\left[r(\om_i)+r(-\om_{i+1})\right]\cdot\left[\vp(\omi)-\vp(\un\om)\right]\right)\right\}=\\
=\Evt\Biggl\{\sum_{i\in\Zb}\biggl(\left[\e{\te_i}\frac{\mu^{(\te_i)}(\om_i-1)}{\mu^{(\te_i)}(\om_i)}+\e{-\te_{i+1}}\frac{\mu^{(\te_{i+1})}(\om_{i+1}+1)}{\mu^{(\te_{i+1})}(\om_{i+1})}\right]\times\\
\times\left[\vp(\omi)-\vp(\un\om)\right]\biggr)\Biggr\}.
\end{multline*}
Under this expectation, the fractions above increase $\om_i$ and decrease $\om_{i+1}$ in the argument of $\vp$ by one, respectively. (This can be seen e.g.\ by changing those variables under $\Evt$.) Remembering what $\omi,\,\omip,\,\omim$ means, this leads us to
\begin{multline*}
\frac{\di}{\di t}\Evt\,\vp(\un\om)=\Evt\Biggl\{\sum_{i\in\Zb}\Bigl[\e{\te_i}\vp(\omipp)-\e{\te_i}\vp(\omip)+\\
+\e{-\te_{i+1}}\vp(\omim)-\e{-\te_{i+1}}\vp(\omipm)\Bigr]\Biggr\}.
\end{multline*}
Now let $l$ and $r$ denote the indices for which the cylinder function $\vp$ depends on only $\om_l,\,\dots,\,\om_r$. Then $\vp(\omip)$ and $\vp(\omim)$ are identical to $\vp(\un\om)$ for $i$'s not between $l$ and $r$ inclusive. Hence
\begin{multline*}
\frac{\di}{\di t}\Evt\,\vp(\un\om)=\Evt\Biggl\{\sum_{i=l-1}^r\Bigl[\e{\te_i}\vp(\omipp)-\e{\te_i}\vp(\omip)+\\
+\e{-\te_{i+1}}\vp(\omim)-\e{-\te_{i+1}}\vp(\omipm)\Bigr]\Biggr\}=\\
=\Evt\Biggl\{\sum_{i=l-1}^{r-1}\left[\e{\te_i}-\e{\te_{i+1}}\right]\vp(\omipp)+\\
+\sum_{i=l+1}^{r+1}\left[\e{-\te_i}-\e{-\te_{i-1}}\right]\vp(\omimm)+\\
+\e{\te_r}\vp(\un\om^{r+1,\,+})-\e{\te_{l-1}}\vp(\un\om^{l-1,\,+})+\e{-\te_l}\vp(\un\om^{l-1,\,-})-\e{-\te_{r+1}}\vp(\un\om^{r+1,\,+})\Biggr\}=\\
=\Evt\Biggl\{\sum_{i=l-1}^{r-1}\left[\e{\te_i}-\e{\te_{i+1}}\right]\vp(\omipp)+\\
+\sum_{i=l+1}^{r+1}\left[\e{-\te_i}-\e{-\te_{i-1}}\right]\vp(\omimm)+\\
+\left[\e{\te_r}-\e{\te_{l-1}}+\e{-\te_l}-\e{-\te_{r+1}}\right]\vp(\un\om)\Biggr\}=\\
=\sum_{i=l-1}^{r-1}\left[\e{\te_i}-\e{\te_{i+1}}\right]\cdot\left[\Evt\,\vp(\omipp)-\Evt\,\vp(\un\om)\right]+\\
+\sum_{i=l+1}^{r+1}\left[\e{-\te_i}-\e{-\te_{i-1}}\right]\cdot\left[\Evt\,\vp(\omimm)-\Evt\,\vp(\un\om)\right].
\end{multline*}
This formula only differs from the right-hand side of \eqref{eq:fo} in indices $i$ where the summands are zero (either because $\te$'s agree or due to the cylinder property of $\vp$).
\end{proof}

\begin{cor}\label{cor:Gibbs}
The product measure $\un\mu^{(\un\te)}$ with $\te_i\equiv\text{constant}$ is time-stationary for the bricklayers' process.
\end{cor}

\noindent
{\bf Remark.} Note that up to this point a computation for the totally asymmetric nearest-neighbor zero range process would look very similar. For these arguments to apply to zero range, one would have begun with $\om_i\in\Zb^+$ in \eqref{eq:Omega}, and put zero to the right-hand side of \eqref{eq:ratafelt}. Then formally we obtain our formulas for zero range by simply neglecting all terms containing $r(-\om)$ and $\e{-\te}$. However, the next lemma requires doubly infinite possibilities for values of $\om_i$, hence, by its essence, it cannot be applied to zero range.

Fix $\un\te$ and
\begin{equation}
\un\te':\,=(\dots,\,\te_{i-1},\,\te',\,\te_{i+1},\,\dots).\label{eq:vesszo}
\end{equation}
We are after the possible values of $\te'$, if any, for which $\Evt\vp(\un\om^{(i,\,\pm)})$ can be written as the expectation of $\vp(\un\om)$ w.r.t.\ the expectation $\Ev^{(\un\te')}$ according to this new parameter-vector $\un\te'$.
\begin{lm}\label{lm:eltol}
With the notation of \eqref{eq:vesszo},
\begin{equation}
\Evt\vp(\un\om^{(i,\,\pm)})=\Ev^{(\un\te')}\vp(\un\om)\label{eq:lep}
\end{equation}
for any cylinder function $\vp$ if and only if the model is the exponential bricklayers' model with parameter $\be=\pm(\te'-\te_i)$ (see \eqref{eq:eblr}). In other words,
\[
\Evt\vp(\un\om^{(i,\,\pm)})=\Ev^{(\un\te^{(i,\,\pm)})}\vp(\un\om)\label{eq:lepp}
\]
holds only for the exponential bricklayers' process, and there is no other set of parameters which can be written in place of $\un\te^{(i,\,\pm)}$ (see \eqref{eq:ipm} for the definition of $\un\te^{(i,\,\pm)}$).
\end{lm}
\begin{proof}
As $\Evt$ and $\Ev^{(\un\te')}$ are product measures, \eqref{eq:lep} is equivalent to
\[
\mu^{(\te_i)}(z\mp1)=\mu^{(\te')}(z)
\]
for all $z\in\Zb$. By the definition \eqref{eq:om} of these measures and property \eqref{eq:ratafelt}, we can expand this equation:
\begin{equation}
\begin{split}
\frac{1}{Z(\te_i)}\cdot\frac{\e{\te_i(z\mp1)}}{r(z\mp1)!}&=\frac{1}{Z(\te')}\cdot\frac{\e{\te'\cdot z}}{r(z)!}\\
r(\pm z)&=\frac{Z(\te_i)}{Z(\te')}\cdot\e{(\te'-\te_i)z}\e{\pm\te_i}.
\end{split}\label{eq:foex}
\end{equation}
From this we see that $r$ cannot be other than exponential of the form \eqref{eq:eblr} with parameter $\be=\pm(\te'-\te_i)$. We need to verify that in this case \eqref{eq:foex} in fact holds. To this order, we assume $\te'=\te_i\pm\be$ and $r(z)=\e{-\be/2+\be z}$. We observe that due to the discrete normal structure, $\widetilde Z(\be,\,m)$ defined in \eqref{eq:emu} is periodic in the variable $m$ with period one. This means
\begin{multline*}
Z(\te_i)=\e{\frac{\te_i^2}{2\be}}\widetilde Z\left(\be,\,\frac{\te_i}{\be}\right)=\e{\frac{(\te'\mp\be)^2}{2\be}}\widetilde Z\left(\be,\,\frac{\te'}{\be}\mp1\right)=\\
=\e{\mp\te'+\frac{\be}{2}}\e{\frac{\te'^2}{2\be}}\widetilde Z\left(\be,\,\frac{\te'}{\be}\right)=\e{\mp\te'+\frac{\be}{2}}Z(\te')=\e{\mp\te_i-\frac{\be}{2}}Z(\te').
\end{multline*}
Plugging this into the right-hand side of \eqref{eq:foex} we have
\[
\frac{Z(\te_i)}{Z(\te')}\cdot\e{(\te'-\te_i)z}\e{\pm\te_i}=\e{-\frac{\be}{2}}\e{(\te'-\te_i)z}=\e{-\frac{\be}{2}+\be(\pm z)},
\]
which agrees to the left-hand side of that equation.
\end{proof}

With this result and definitions \eqref{eq:ipm}, \eqref{eq:ipmo} in mind, \eqref{eq:fo} can be written in the form
\begin{multline*}
\frac{\di}{\di t}\Evt\,\vp(\un\om)=\sum_{i:\,\te_i\ne\te_{i+1}}\left(\e{\te_i}-\e{\te_{i+1}}\right)\cdot\left(\Evtipp\,\vp(\un\om)-\Evt\,\vp(\un\om)\right)+\\
+\sum_{i:\,\te_{i-1}\ne\te_i}\left(\e{-\te_i}-\e{-\te_{i-1}}\right)\cdot\left(\Evtimm\,\vp(\un\om)-\Evt\,\vp(\un\om)\right)
\end{multline*}
for the exponential bricklayers' process, and no other bricklayers' (nor totally asymmetric nearest-neighbor zero range) processes. Hence we have arrived to the proof of Theorem \ref{tm:fo}.

\hop
{\bf Remark.} For other than exponential rates, Lemma \ref{lm:fo} tells us that $\un\mu^{(\un\te)}$ evolves into a product of our original marginals $\mu^\te(\cdot)$ and of these marginals shifted upwards or downwards by one: $\mu^\te(\cdot\pm1)$. However, these shifted measures do not correspond to the original marginals with shifted parameters as Lemma \ref{lm:eltol} does not apply for other than exponential models. Therefore the class of product measures of marginals $\mu^\te$ is not closed under the time-evolution in this case. As we do not know how the product of shifted versions (in the sense above) of $\mu^\te$ evolves, we do not have any positive result for the non-exponential case.

\begin{proof}[Proof of lemma \ref{lm:w}] For the purposes of this proof, we enumerate all $n$ random walkers. At each moment, we label them from $1$ to $n$ in such a way that a walker right to another must have larger label than the other one. We translate their jump rates according to this labeling: while the $m$-th walker is alone at his site, she has rate $\e{\te_\text{left}-(m-1)\be}-\e{\te_\text{left}-m\be}$ to jump to the right, and rate $\e{-\left(\te_\text{left}-m\be\right)}-\e{-\left(\te_\text{left}-(m-1)\be\right)}$ to jump to the left, respectively. If more of them meet at a site then the lowest labeled of them has the sum of these left-rates as rate to jump left with, and the largest labeled of them has the sum of the right-rates to jump to the right with, respectively (see the remark after \eqref{eq:sokLR}). All others have zero rate to jump due to our labeling method. These rates are in accordance with the rates defined in Section \ref{sc:muso}.

Define $\eta_m$ as the distance between the $m$-th and the $m+1$-th walker ($m=1\dots n-1$). Then $\eta_m\ge0$, and, whenever it is positive, it decreases by one with rate larger than or equal to
\begin{multline*}
\e{\te_\text{left}-(m-1)\be}-\e{\te_\text{left}-m\be}+\e{-\left(\te_\text{left}-(m+1)\be\right)}-\e{-\left(\te_\text{left}-m\be\right)}=\\
=\e{\be}\left(1-\e{-\be}\right)\left(\e{\te_\text{left}-\be m}+\e{-\te_\text{left}+\be m}\right).
\end{multline*}
This would be the rate if the $m$-th and $m+1$-th walkers were alone at their sites. If they are not alone, then by the labeling method, $m$ is the first one to jump to the right and $m+1$ is the first to jump to the left, hence $\eta_m$ decreases even faster. However, when $\eta_m$ is positive, its maximum rate to increase is
\begin{multline*}
\e{-\left(\te_\text{left}-m\be\right)}-\e{-\left(\te_\text{left}-(m-1)\be\right)}+\e{\te_\text{left}-m\be}-\e{\te_\text{left}-(m+1)\be}=\\
=\left(1-\e{-\be}\right)\left(\e{\te_\text{left}-\be m}+\e{-\te_\text{left}+\be m}\right).
\end{multline*}
This would be rate for $\eta_m$ to grow if the walkers $m$ and $m+1$ were alone at their sites. If they are not, then by the labeling method $m$ is the last one to jump left and $m+1$ is the last one to jump right, so this rate is even smaller.

Thus we see that $\eta_m$ is a process on $\Zb^+$ which increases with smaller rates than decreases (except when it is zero). Hence its stationary distribution is stochastically bounded by a geometric one (disturbed at zero) and so all its moments are finite. As the width of the shock is just
\[
\sum_{m=1}^{n-1}\eta_m,
\]
we see that it also has finite moments.
\end{proof}
\begin{proof}[Proof of lemma \ref{lm:cm}] Jumps of walkers from different positions happen independently, with rates
\[
\begin{array}{rl}
\e{\te_\text{L}}-\e{\te_\text{R}}&\text{to the right, and}\\
\e{-\te_\text{R}}-\e{-\te_\text{L}}&\text{to the left}
\end{array}
\]
with $\te_\text{L}$ and $\te_\text{R}$ being the local parameter values to the left and to the right of the position in question. The rate with which \emph{any} of our $n$ walkers jumps is therefore the sum of the previous rates, and is
\begin{equation}
\begin{array}{rl}
\e{\te_\text{left}}-\e{\te_\text{right}}&\text{to the right and}\\
\e{-\te_\text{right}}-\e{-\te_\text{left}}&\text{to the left},
\end{array}\label{eq:cmr}
\end{equation}
where $\te_\text{left}$ ($\te_\text{right}$) is the value of the parameter left to the leftmost (right to the rightmost, respectively) walker.

$X(t)$ defined in \eqref{eq:cm} is a rescaled sum of the walkers' position, therefore it increases (decreases) by $1/n$ if and only if \emph{any} of our walkers jumps to the right (left, respectively). This happens according to a Poisson-process, with rates shown in \eqref{eq:cmr}.
\end{proof}

\nemfej{Acknowledgment}

The author wishes to thank G\"unter Sch\"utz for drawing his attention to the possibility of this phenomenon, Timo Sepp\"al\"ainen and Firas Rassoul-Agha for illuminating conversations about the subject. We also thank the referees for careful reading of the previous version of the paper, and pointing out the parts where clarity could be improved.

\bibliography{eredeti}
\bibliographystyle{plain}
\bigskip\bigskip\bigskip\bigskip
\begin{flushright}
{\small\begin{tabular}{l}
M\'arton Bal\'azs\\
{\tt balazs@math.wisc.edu}\\
Mathematics Department,\\
University of Wisconsin-Madison\\
480 Lincoln Dr\\
Madison WI 53706-1388\\
USA
\end{tabular}}
\end{flushright}
\end{document}